\documentclass[a4paper]{article}
\usepackage{amsmath,amssymb,amsfonts}

\newtheorem{definition}{Definition}[section]

\title{On power sets}

\author{Jailton C. Ferreira}

\date{ }
\begin{document}
\maketitle
\pagenumbering{arabic}
\begin{abstract}
This work presents theorems which state (i) $Z$ is a proper subset
for any bijection $f$ between $A$ and $Z$, where $Z \subseteq
P(A)$, $A$ is a non-finite set and $|Z|=|A|$, and (ii) being $Z$ a
proper subset of $P(A)$ nothing affirms or denies that
$|P(A)|>|A|$. Russell's paradox is examined and it is shown that
the set of all the ordinary sets does not exist. A mistake in
Cantor's proof on cardinality of power sets is shown
\end{abstract}
\section{Introduction}

\hspace{22pt} We cannot decide, applying the diagonal argument, if
it is true that the cardinality of the set of real numbers of
interval [0,1] is larger than the set of natural numbers $N$
~\cite{Ferreira}. The diagonalization argument appears clearly in
proofs of $|P(N)|>|N|$ in which the characteristic functions of
the subsets of $N$ are used. Consequently, if the results of
~\cite{Ferreira} are true, we also cannot decide, applying the
diagonal argument, if $|P(N)|>|N|$ is true. In section \textbf{2}
two theorems state (i) $Z$ is a proper subset for any bijection
$f:A \rightarrow Z$ where $Z \subseteq P(A)$, A is a non-finite
set and $|Z|=|A|$, and (ii) being $Z$ a proper subset of $P(A)$
nothing affirms or denies that $|P(A)|>|A|$.

\hspace{22pt} Section \textbf{3} examines Russell's paradox and
shows that the set of all ordinary sets does not exist. This
section is included to help the argumentation of the next section.

\hspace{22pt} Section \textbf{4} shows a mistake in Cantor's proof
on cardinality of power sets.

\section{Theorems on power sets}

\subsection{Cantor's theorem}

\emph{Theorem}. The set $P(A)$ of all subsets of a set $A$ has a
larger cardinality than $A$.

\emph{Proof}. Suppose they have the same number of elements.

Let $b:A \rightarrow P(A)$ be a bijection between $A$ and $P(A)$.

(1) Let  $T=\{x$ in $A | x \notin b(x)\}$.

Since $T$ is a subset of $P(A)$ and $b$ is onto,

(2)  $T=b(t)$  for some $t$.

Thus  $t$ is in $b(t)$  iff (by 2)  $t$ is in $T$ iff (by 1)  $t$
is not in $b(t)$.

This is a contradiction. Since we cannot do an one-to-one
correspondence between $A$ and $P(A)$, and since $P(A)$ cannot be
smaller than $A$, the only possible conclusion is that $P(A)$ must
be larger than $A$.

\subsection{Theorem 1}
\begin{em}
\hspace{22pt} Let  $f:A \rightarrow Z$ be a bijection between $A$
and $Z$, where $A$ is a non-finite set, $Z$ is a subset of $P(A)$
and

\begin{equation}
|Z|=|A|
\end{equation}

The set $Z$ is a proper subset of $P(A)$.
\end{em}

\emph{Proof}.

\hspace{22pt} The set
\begin{gather}
\text{$Y=\{x $ in $A | x \notin f(x)\}$}
\end{gather}

must exist by the axiom of specification, $P=\{x$ in $Q
|Prop(x)\}$. Consider some element $a$ of $A$ that does not belong
to the set $f(a)$, that is, $Y$ is not empty.

\hspace{22pt} Since $Y \in P(A)$ and $f$ is bijective, there is no
$c$ belonging to $A$ such that

\begin{equation}
f(c) = Y
\end{equation}

The hypothetic element $c$ does not belong to $A$ because it
cannot belong to $Y$. Therefore

\begin{equation}
Z \subset P(A)
\end{equation}

that is, the set $Z$ is a proper subset of $P(A)$.\\

\subsection{Theorem 2}

\begin{em}
\hspace{22pt} Let $Z$ be a subset of $P(A)$ such that $f:A
\rightarrow Z$ is a bijection between $A$ and $Z$, $A$ is a
non-finite set and $|Z|=|A|$. Being $Z$ a proper subset of $P(A)$
nothing affirms or denies that $|P(A)|=|A|$.
\end{em}\\

\emph{Proof}.

\hspace{22pt} Let $f:A \rightarrow P(A)$ be a injection between
$A$ and $P(A)$. Let us consider an element $Y$ of $P(A)$ such that
\begin{equation}
Y \neq \varnothing
\end{equation}
\begin{equation}
Y \cap f(A) = \varnothing
\end{equation}
and
\begin{gather}
\text{$Y=\{x $ in $A | x \notin f(x)\}$}
\end{gather}
The bijective function $f:A \rightarrow f(A)$ satisfies the
Theorem 1. However we do not know if $|P(A)|>|A|$ is true or false
because the cardinality of
\begin{equation}
D = P(A) - (f(A) \cup Y)
\end{equation}
is ignored. Since
\begin{equation}
|f(A) \cup Y| = |A|
\end{equation}
we obtain
\begin{equation}
|P(A)| = |D| + |A|
\end{equation}

\section{The paradox of the sets}
\begin{quotation}
Let ``$x$ is a set that is not member of $x$" and $O$ the set that
it determines, set $O$ will be member of $O$ if and only if it
satisfies the function that determines $O$, that is, if and only
if $O$ is not member of $O$, what is a contradiction.
\end{quotation}

This contradiction is known as the paradox of the sets or, more
generally, as Russell's paradox. To examine the contradiction we
will use the usual definitions:

\begin{definition}
A set is ordinary if it is not member of itself.
\end{definition}

\begin{definition}
An extraordinary set is a set that belongs to itself.
\end{definition}

\hspace{22pt} Let us consider the following work hypothesis: given
any property, there is a set of all things that have this
property. Let the property of being an ordinary set. Applying the
work hypothesis, we have a set of all ordinary sets; let us
denominate $H$ the set. Let us now consider the following
argument:

\begin{enumerate}
\item[(a)] If $H$ does not belong to $H$, then $H$ is an ordinary set, but
then $H$ does not contain all the ordinary sets.
\item [(b)] If $H$ belongs to $H$, then $H$ is an extraordinary set, but then $H$ contains
a set that is not ordinary.
\end{enumerate}

There are two necessary conditions for $H$: (i) to contain all
ordinary sets; and (ii) not to contain any extraordinary set. From
(a) and (b), we see that the conditions are not satisfied: the set
$H$ cannot belong to $H$ to exclude extraordinary set in $H$, and
$H$ cannot exclude $H$ to contains all the ordinary sets of $U$.
If $H$ exists, then
\begin{equation}
H \supset S
\end{equation}

where $S$ is an ordinary set. A no empty set cannot satisfy both
conditions (i) and (ii). Therefore $H$ does not exist.

\section{On the sets T and Y}
\hspace{22pt} Let us consider a no empty proper subset of $A$,
denoted by $B$, such that their elements satisfy ``$x$ is not in
$b(x)$" where $b:B \rightarrow W$ and $W$ is a subset of $P(A)$.
If $T$ of Cantor's theorem exists, then
\begin{equation}
T \supset B
\end{equation}
and, therefore, $T$ is not empty.

\hspace{22pt} The set $T$ is the set of the elements $x$ of $A$
such that ``$x$ is not in $b(x)$" (the elements of $T$ are the
elements of $A$ which separately satisfy the property ``$x$ is not
in $b(x)$") and \\

\hspace{80pt}\text{``$T$ contains the element $a$ satisfying
$b(a)=T$ and}
\begin{gather}\label{beta}
\text{$T$ does not contain the element $a$ of $A$ satisfying
$b(a)=T$"}
\end{gather}

A no empty set cannot satisfy the condition (\ref{beta}).
Therefore $T$ does not exist.

\hspace{22pt} As $T$ does not exist, the contradiction of Cantor's
theorem, section \textbf{2.1}, does not exist and we cannot
conclude
\begin{equation}
|P(A)|>|A|
\end{equation}

\hspace{22pt} In the case of Theorem 1 there is not a paradox, but
an impossibility ``The hypothetic element c does not belong to $A$
because it cannot belong to $Y$".

\end{document}